\documentclass{article}
\usepackage{amsmath,amssymb,amsthm}

\begin{document}
\title{\bf Finite simple groups with some abelian Sylow subgroups
\thanks{Project supported by NSFC (No.11201133 and 11471131).}}

\author{{Rulin Shen$^{1}$~~~~  Yuanyang Zhou$^2$}
\\  {\small
  $^1$Department of Mathematics,
  Hubei  University for Nationalities, }\\ {\small Enshi, Hubei Province,
   445000, P. R. China}\\{\small
  $^2$Department of Mathematics,
  Central China Normal University, }\\ {\small Wuhan, Hubei Province,
   430079, P. R. China} }

\date{}\maketitle

\insert\footins{\footnotesize{\it Email addresses}:\
rshen@hbmy.edu.cn (Rulin Shen) }

\begin{abstract}
In this paper, we classify the finite simple groups with an
 abelian Sylow  subgroup.\medskip

\medskip

\end{abstract}
\section{Introduction}

Sylow subgroups are very important subgroups of finite groups. Some
well-known theorems on the structure of  groups link their Sylow
subgroups. In particular,   groups with some abelian Sylow subgroups
have been widely researched. At the same time, some problems related
to abelian Sylow subgroups are put forward. Since the condition of
abelian Sylow subgroups inherits to subgroups and quotient groups,
some problem will reduce to  simple having some abelian Sylow
subgroups.

In 1969, J. H. Walter   classified   finite non-abelian simple
groups with abelian Sylow 2-subgroups (ref. \cite{Wa}), which are
$L_2(q)$, where $q=2^f (f\ge 2)$ or $q\equiv 3,5(\bmod~8)$, ${ J}_1$
or $^2G_2(q)$, where $q=3^{2m+1}$ and $m\ge 1$. Also the structure
of Sylow subgroups is studied by A. J. Weir(ref. \cite{Weir}). He
proved that the Sylow $r$-subgroups of the classical groups over
$\rm{GF}$$(q)$ with $q$ prime to $r$ and $r$ odd are expressible as
direct products of groups defined inductively by $G_i := G_i\wr C_p
$, where $G_0$ is an abelian $p$-subgroup of the classical group and
$C_p$ a cyclic group of order $p$.

In this paper, we use Artin invariant of the cyclotomic
factorisation(ref. \cite{KLST}) to classify the finite simple groups
with an
 abelian Sylow  subgroup.\\

 \noindent\textbf{Theorem.} Let $S$ be a finite non-abelian simple
 group. Suppose that $S$ has an abelian Sylow $r$-subgroup and $r$ an odd prime.  Then
 $S$ is one of the following groups:

 (1) an alternating group $A_n$ with $n<r^2$,

 (2) a linear simple group ${\rm PSL}_2(q)$ with $r|q$,

 (3) ${\rm PSL}_3(q)$,  $q\equiv 4,7(\bmod~9)$  with $r=3$,

 (4) ${\rm PSU}_3(q^2)$,  $2<q\equiv 2,5(\bmod~9)$ with $r=3$,

  (5) a simple group  of Lie type
over the Galois field $\rm{GF}$$(q)$ and $ {e}_L(mr_m)=0$ except the
above groups of (3) and (4), where $r=r_m$ is a primitive prime of
$q^m-1$, the function $e_L(x)$ is defined in Table 1 and Table 2,

(6) a sporadic simple groups listed in Table 4. \\

From above theorem, we can get the following result.\\

\noindent\textbf{Corollary.} Let $S$ be a finite  simple
 group. Suppose that $S$ has an abelian Sylow $r$-subgroup $R$. Then
 $R$ is isomorphic to some direct products of copies of a cyclic $r$-group.\\

 Afterwards we   consider the classification of finite simple groups having an elementary abelian Sylow subgroup. In the following, all considered groups are finite. Recall that
 a \emph{section} of a group $G$  is a quotient of a subgroup of $G$. Denote by $G^k$ the group of $G\times \cdots \times G$ with $k$ times. The notations
 $G.K$ and $G:K$   mean an extension and a split  extension  of
 the group $G$ by $K$, respectively.  We denote $p^s\parallel n$ if $p^s | n$ but $p^{s+1}\nmid n$ for a prime $p$ and natural number $n$.\medskip

\section{Alternating Groups}

\noindent\textbf{Theorem 2.1.}  Let $r$ be an odd prime.  Then
 the alternating group $A_n$  has a Sylow abelian $r$-subgroup   if and
only if $n<r^2$. \medskip

\textbf{Proof.}  Since $r$ is an odd prime, Sylow $r$-subgroups of
$A_n$ are ones of the symmetric group ${\rm Sym}(n)$.  We use the
well known structure of Sylow subgroups of symmetric groups. Now let
$r^s$ be the largest number such that $r^s\leq n$. If
$[\frac{n}{r^s}]\geq 1$ and $s\geq 2$, then there exists a subgroup
$C_r\wr C_r$, which is not abelian. So Sylow $r$-subgroups of $A_n$
are not abelian. Next we assume that $s=1$ and $[\frac{n}{r}]\leq
r-1$. Then Sylow $r$-subgroups of $A_n$ is $C_r^{[\frac{n}{r}]}$,
which is an elementary abelian group.\hfill$\blacksquare$\\

Next we give a result on the alternating groups having a Sylow
$r$-subgroup isomorphic to $C_r\times C_r$.\\

\noindent\textbf{Corollary  2.2.} Let $r$ be an odd prime.  Then the
alternating group $A_n$ $(n\geq 5)$ has a Sylow $r$-subgroup $C_r^2$
if and only if $2r\leq n< 3r.$\medskip

\section{Simple Groups of Lie Type}
In this section, we will discuss simple groups of Lie type over the
Galois field {\rm GF}$(q)$.  For any natural number $m$ let
$\Phi_m(x)$ denote the $m$th cyclotomic polynomial. So the degree of
$\Phi_m(x)$ is the Euler function $\phi(m)$ and $x^m-1=
\prod_{d|m}\Phi_d(x)$ (ref. \cite[Page 207]{R}). Recall that a
primitive prime divisor (or Zsigmondy prime) $r_m$ of $q^m-1$ is a
prime such that $r_m | (q^m-1)$ but $r_m \nmid\ (q^i-1)$ for $1\leq
i\leq m-1$. A well-known theorem by Zsigmondy asserts that primitive
primes of $q^m-1$ exist except if $(q,m)=(2,6)$ or $m=2$ and
$q=2^k-1$ (ref. \cite{Z}).
\\

\noindent\textbf{Lemma 3.1.} Let $q$ and $m$ be   natural numbers.
Then $\Phi_1(q)\equiv -1(\bmod~ q)$ and  $\Phi_m(q)\equiv 1(\bmod~
q)$ for $m\ge 2$.\medskip

\textbf{Proof.} Clearly the result is true for $m=1$ and 2. Next we
prove the remaining case by induction on $m$. Since
$q^m-1=\Phi_1(q)\Phi_{m}(q)\prod_{d|m, d\neq 1,m}\Phi_d(q)$,   we
have $q^m-1\equiv \Phi_1(q)\Phi_{m}(q)\equiv -\Phi_{m}(q)(\bmod~ q)$
by induction, and hence
$\Phi_{m}(q)\equiv 1-q^m\equiv 1(\bmod~ q)$ whenever $m\geq 2$. \hfill$\blacksquare$\\

\noindent\textbf{Lemma 3.2} (\cite[Page 207]{R}).  Let $q, m$ be
positive integers with $q\ge 2$ and $r_m$ a primitive prime of
$q^m-1$. Then $r_m$ divides $\Phi_m(q)$.\\

\noindent\textbf{Lemma 3.3} (\cite[Lemma 5]{M}). Let $r_m$ be an odd
primitive prime of $q^m-1$. Then $r_m| \Phi_n(q)$ if and only if
$n=mr_m^j$ for some $j\geq 0$. \\

\noindent\textbf{Lemma 3.4.} Let $r_m$ be an odd primitive prime of
$q^m-1$ and $j\geq 1$. Then $r_2^j\parallel\Phi_{2r_2^j}(q)$ and
$r_m|| \Phi_{mr_m^j}(q)$ for $m\neq 2$.\medskip

\textbf{Proof.} The proof of the case $m\geq 3$ is given by W. Feit
(ref. \cite[Lemma 2.1]{F}) or E. Artin (ref.  \cite[Lemma 1]{A}).
Next we prove the cases   $m=1$ and 2.

First, we consider the case $m=1$. Now we set $q-1\equiv
k(\bmod~r_1^2)$ with $0\leq k\le r_1^2-1$. Since $r_1|(q-1)$, we
have $r_1|k$, and let $k=lr_1$. Then $q\equiv lr_1+1(\bmod~r_1^2)$
where $0\le l\le r_1-1$.
 Moreover
$$\Phi_{r_1^j}(q)=\frac{q^{r_1^j}-1}{q^{r_1^{j-1}}-1}=q^{r_1^j-r_1^{j-1}}+q^{r_1^j-2r_1^{j-1}}+\cdots + q^{r_1^j-(r_1-1)r_1^{j-1}}+1,$$
so that $\Phi_{r_1^j}(q)\equiv r_1(\bmod~r_1^2)$, and then
$r_1\parallel \Phi_{r_1^j}(q)$.

We next discuss the case $m=2$. Similarly, we let $q+1\equiv
k(\bmod~r_2^{j+1})$ with $0\leq k\le r_2^{j+1}-1$. Since
$r_2|(q+1)$, we have $r_2|k$, and let $k=lr_2$. Then $q\equiv
lr_2-1(\bmod~r_2^{j+1})$ where $0\leq l\leq r_2^j-1$.  Since
$\Phi_{2r_2^j}(q)=\sum_{i=0}^{r_2^j-1}(-1)^iq^i,$
 it follows that
$$\Phi_{2r_2^j}(q)\equiv\sum_{i=0}^{r_2^j-1}(-1)^i(lr_2-1)^i=\frac{(lr_2-1)^{r_2^j}+1}{lr_2}=\sum_{i=1}^{r_2^j}(-1)^{i-1}(^{r_2^j}_i)(lr_2)^{i-1}$$
$$~~~~~~\equiv r_2^j(\bmod~r_2^{j+1}),$$
and thus $r_2^j\parallel\Phi_{2r_2^j}(q)$.\hfill$\blacksquare$\\

In the following, let $L(q)$ be a simple group of Lie type of
charactersic $p$. The order of $L(q)$ has the cyclotomic
factorization in terms of $q$:
$$|L(q)|=\frac{1}{d}q^h\prod_m\Phi_m(q)^{e_L(m)}\eqno (1)$$
where $d$ is the denominator and $h$ the exponent given for $L(q)$
in Table 1 and Table 2, where $\Phi_m(q)$ is the cyclotomic
polynomial for the primitive $m$th roots of unity, and where the
$e_L(m)$ are the exponents deducible from the factors. The function
$e_L(x)$ is defined in Table 1 and Table 2 (ref. \cite[Table
C.1,C.2]{KLST}).\medskip

\begin{small} ~~~~~~Table 1: {\bf~Cyclotomic factorisation: classical groups of Lie
type}\smallskip

\begin{tabular}{|l|l|l|l|} \hline
\quad $L$ &$ d$ &$ h$ & $ e_L(x)$ where $x$ is an integer \\[2pt]
\hline
 ${\rm PSL}_{n}(q)$ & $(n,q-1)$   & $n(n-1)/2$& $n-1~~~$if $x=1$\\[2pt]
 &&&$[\frac{n}{x}]~~~$ if $x>1$ \\[2pt]\hline
&&&$[\frac{n}{{\rm lcm}(2,x)}]$~ if $x\not\equiv 2(\bmod~4)$\\[2pt]
${\rm PSU}_n(q^2)$&$(n,q+1)$&$n(n-1)/2$& $n-1$ ~ if $x=2$\\[2pt]
&&&$[\frac{2n}{x}]$ ~ if $2<x \equiv 2(\bmod~4)$\\[2pt]\hline
 ${\rm PSp}_{2n}(q), {\rm P\Omega}_{2n+1}(q)$ &$(2,q-1)$ &$n^2$&
$[\frac{2n}{{\rm lcm}(2,x)}]$\\[4pt]\hline
${\rm P\Omega}^+_{2n}(q)$&$(4,q^n-1)$&$n(n-1)$&$[\frac{2n}{{\rm
lcm}(2,x)}]$ ~if
$x|n$ or $x\nmid 2n$\\[4pt]
&&& $\frac{2n}{x}-1$~ if $x\nmid n$ and $x| 2n$
\\[2pt]\hline
${\rm P\Omega}^-_{2n}(q)$&$(4,q^n+1)$&$n(n-1)$&$[\frac{2n}{{\rm
lcm}(2,x)}]-1$ ~if
$x|n$ \\[2pt]
&&& $[\frac{2n}{{\rm lcm}(2,x)}]$ ~if $x\nmid n$
\\[2pt]\hline
\end{tabular}\medskip\medskip

  ~~Table 2: {\bf~Cyclotomic factorisation: exceptional groups of Lie
type}\smallskip \\
~\begin{tabular}{|l|l|l|l|l|l|l|l|l|l|l|} \hline
 $m$& $^2{\rm B}_2$ &$^3{\rm D}_4$ &$ {\rm G}_2$ & $ ^2{\rm G}_2$
&${\rm F}_4$&$^2{\rm F}_4$&${\rm E}_6$&$^2{\rm E}_6$&${\rm
E}_7$&${\rm E}_8$
\\\hline
1&1&2&2&1&4&2&6&4&7&8\\\hline 2&&2&2&1&4&2&4&6&7&8\\\hline
3&&2&1&&2&&3&2&3&4\\\hline 4&1&&&&2&2&2&2&2&4\\\hline
5&&&&&&&1&&1&2\\\hline 6&&2&1&1&2&1&2&3&3&4\\\hline
7&&&&&&&&&1&1\\\hline 8&&&&&1&&1&1&1&2\\\hline
9&&&&&&&1&&1&1\\\hline 10&&&&&&&&1&1&2\\\hline
12&&1&&&1&1&1&1&1&2\\\hline 14&&&&&&&&&1&1\\\hline
15&&&&&&&&&&1\\\hline 18&&&&&&&&1&1&1\\\hline 20&&&&&&&&&&1\\\hline
24&&&&&&&&&&1\\\hline 30&&&&&&&&&&1\\\hline
$d$&1&1&1&1&1&1&$(3,q-1)$&$(3,q+1)$&$(2,q-1)$&1\\\hline
$h$&2&12&6&3&24&12&36&36&63&120\\\hline
\end{tabular}\end{small}

Note that for Table 2, the numbers in the row headings are those
integers $m$ for which a cyclotomic polynomial $\Phi_m(q)$ enters
into the cyclotomic factorisation in terms of $q$ of one of the
exceptional groups $L(q)$ of Lie type according to the formula (1).
The types $L$ are the headings, and the entry in the $m$-row and the
$L$-column is $e_L(m)$ if non-zero and blank otherwise.

 Moreover, let $R_m$ be the set of all primitive
primes of $q^m-1$. Denote by $Z_m(q)=\prod_{r_m\in R_m}
|\Phi_m(q)|_{r_m}$, where $|\Phi_m(q)|_{r_m}$ means the $r_m$'s
maximum prime power divisor of $|\Phi_m(q)|$. Then we have
$$|L(q)|=\frac{1}{d}q^h\prod_m\prod_{r_m\in R_m} Z_m(q)^{e_L(m)}r_m^{\widetilde{e}_L(m)}\eqno(2)$$
where  $\widetilde{e}_L(m)$ is the number of divisor $r_m$ of all
possible
$\Phi_{m'}(q)$ such that $m'>m$. Note that we define  $Z_6(2)=1$ and $Z_{2}(r)=1,$ whenever $r$ is a Mersenne  prime. \\

\noindent\textbf{Lemma 3.5.}  Keep the above notations. Suppose that
$r_m$ is an odd prime. Then
\[\widetilde{e}_L(m)=\left\{
\begin{array}{ll}\sum_{j=1}^{\infty}e_L(mr_m^j) & \mbox{ if $ m=1$ or $m\geq
3$,
 }\\[6pt]
\sum_{j=1}^{\infty}je_L(mr_m^j)  & \mbox{ if $m=2$}
\end{array}
\right.\]
 for the simple groups of Lie type, where
 the function $e_L(x)$ is
given by Table 1 and Table 2.
\medskip

\textbf{Proof.}  We consider the general term $\Phi_m(q)^{e_L(m)}$
to give a proof. By Lemma 3.3 and Lemma 3.4, it follows that the
$r_m$-part of $\Phi_{mr_m^j}(q)$ is $r_m$ whenever $m\ne 2$, and
$r_m^j$ for $m=2$. So we have $\Phi_{mr_m^j}(q)^{e_L(mr_m^j)}$
$=(r_mZ_{mr_m^j}(q))^{e_L(mr_m^j)}$, if $m\ne 2$, and equals
$(r_m^jZ_{mr_m^j}(q))^{e_L(mr_m^j)}$ if $m=2$. Thus
$\widetilde{e}_L(m)=\sum_{j=1}^{\infty}e_L(mr_m^j)$  for $m\neq 2$,
and $\sum_{j=1}^{\infty}je_L(mr_m^j)$ for $m=2$.
\hfill$\blacksquare$\\

We compute directly in term of Lemma 3.1 and the formulas of Lemma
3.5, the
following result can be obtained.\\

 \noindent\textbf{Lemma 3.6.} Keep
the above notations. Suppose that $r_m$ is an odd primitive  prime
of $q^m-1$. Then $\widetilde{e}_L(m)=0$ for the exceptional simple
groups of Lie type except the  cases listed in Table 3.\\

\begin{small}
~~~~Table 3: {\bf~The function $\widetilde{e}(m)$ for exceptional
groups of Lie type}\smallskip

 ~~~~~~~~~~~~~~\begin{tabular}{|l|l|l|l|l|l|l|l|l|l|}
\hline $r_m$ &$^3{\rm D}_4$ &$ {\rm G}_2$ &${\rm F}_4$&$^2{\rm
F}_4$&${\rm E}_6$&$^2{\rm E}_6$&${\rm E}_7$&${\rm E}_8$
\\\hline
$r_1=3$&2&1&2&&4&2&4&5\\\hline $r_1=5$&&&&&1&&1&2\\\hline
$r_1=7$&&&&&&&1&1\\\hline $r_2=3$&2& &2&1&2&4&4&5\\\hline
$r_2=5$&&&&&&1&1&2\\\hline $r_2=7$&&&&&&&1&1\\\hline
$r_4=5$&&&&&&&&1\\\hline
\end{tabular}
\end{small}
\\

Next we classify  the simple groups of Lie type of   characteristic
$p$ which have an abelian Sylow $p$-subgroup. \\

\noindent\textbf{Theorem 3.7.} Let $L$ be a simple group of Lie type
over the finite field {\rm GF}$(q)$ of characteristic $p$. Then $S$
has an abelian Sylow $p$-group if and only if $L$ is isomorphic to
{\rm PSL}$_2(q)$.
\medskip

\textbf{Proof.} It is known that the Sylow $p$-subgroup of $L$ is a
maximal unipotent subgroup. We use the notations of \cite{V}, let
$a(\Phi,p)$ be the maximum of orders of abelian $p$-subsets of the
root system $\Phi^+$ (also ref. \cite{VV}).  By hypothesis, the
Sylow $p$-subgroup is abelian, so we have  $q^h\leq q^{a(\Phi,p)}$
for classical simple groups (where $a(\Phi,p)$ is listed in
\cite{VV}), and $q^h\leq q^{a(U)}$ for exceptional simple groups
(where $a(U)$ in Table 4 of \cite{V}). Therefore, $L$ has an abelian
Sylow $p$-group  if and only if $L$ is isomorphic to ${\rm
PSL}_2(q)$.\hfill$\blacksquare$\medskip

Next we give a necessary and sufficient condition for a simple group
of Lie type to have an abelian Sylow $r$-subgroups and $r$ is prime
to the characteristic.\medskip

\noindent\textbf{Theorem 3.8.} Let $L$ be a simple group of Lie type
over the Galois field {\rm GF}$(q)$ and $r_m$ an odd primitive prime
of $q^m-1$. Then $L$ has an abelian Sylow $r_m$-subgroup if and only
if $ {e}_L(mr_m)=0$ except the following two cases:

(1) $r_1=3$,  $L\cong {\rm PSL}_3(q)$ and $q\equiv 4,7(\bmod~9)$;

(2) $r_2=3$,  $L\cong {\rm PSU}_3(q^2)$ and $2<q\equiv
2,5(\bmod~9)$.
\medskip

\textbf{Proof.} We use R. W. Carter's results on maximal tori of
simple groups of Lie type in his papers (ref. Propositions 7-10 of
\cite{C1} and \cite{C2}) to prove. We divide into  several
cases.\medskip

\emph{Case $1$: $L\cong {\rm PSL}_n(q)$ with $n\geq 2$}. Then every
maximal torus $T$ of ${\rm PSL}_n(q)$ has  order
$\frac{1}{(n,q-1)(q-1)}(q^{n_1}-1)(q^{n_2}-1)\cdots (q^{n_k}-1)$ for
appropriate partition $n=n_1+n_2+\cdots +n_k$.  We first suppose
that $m\geq 2$. Then there exists a maximal torus has the order
$\frac{1}{(n,q-1)(q-1)}(q^{m}-1)^{[\frac{n}{m}]}$ in terms of the
partition $n=m+m+\cdots+m+l$, where $0\le l<m$, so that the largest
order of abelian $r_m$-subgroups is the $r_m$-part of
$(q^{m}-1)^{[\frac{n}{m}]}$. Thus Sylow $r_m$-subgroups are abelian
if and only if $\widetilde{e}_{L}(m)=0$, which is equivalent to the
condition $e_{L}(mr_m)=0$ in terms of Lemma 3.5. Next we assume
$m=1$, that is $r_1|(q-1)$. We split it into two cases, $r_1$ does
not divide $n$, and then $r_1$ divides $n$ as well.

Subcase 1:  $r_1$ does not divide  $n$.
 Then there
exists a maximal torus of order $\frac{1}{(n,q-1)}(q-1)^{n-1}$,
whose $r_1$-part is the maximum order of the abelian $r_1$-subgroup.
Thus Sylow $r_1$-subgroups are abelian if and only if
$\widetilde{e}_{L}(1)=0$, that is the condition $e_{L}(r_1)=0$.

Subcase 2:  $r_1$   divides  $n$. Let $r_1^s\parallel\gcd(n, q-1)$
and $r_1^t\parallel (q-1)$.  We first consider the group ${\rm
GL}_n(q)$. The Sylow $r_1$-subgroup of  ${\rm GL}_n(q)$ is
isomorphic to one of the group $M$ of monomial matrices which have
one nonzero entry in each row and column, so they are the product of
a permutation matrix and a diagonal matrix, that is $M\cong D:{\rm
Sym}(n)$, where $D$ is the group of the set diagonal matrices. Thus
a Sylow $r_1$-subgroup of ${\rm SL}_n(q)$ is a Sylow subgroup of
$M_1=D_1:{\rm Sym}(n)$, where $D_1$ is the set all diagonal matrices
of determinant  1. Therefore, a Sylow $r_1$-subgroup of $M_1/Z({\rm
SL}_n(q))$ is a Sylow $r_1$-subgroup of ${\rm PSL}_n(q)$. Moreover,
we note that $M_1/Z({\rm SL}_n(q))\cong D_1/Z({\rm SL}_n(q)):{\rm
Sym}(n)$. Obviously, $D_1/Z({\rm SL}_n(q))$ is an abelian group of
order $\frac{(q-1)^{n-1}}{(n,q-1)}$, so that the Sylow
$r_1$-subgroup of ${\rm PSL}_n(q)$ is abelian if and only if   the
$r_1$-parts of $|D_1/Z({\rm SL}_n(q)):{\rm Sym}(n)|$ and $n!$ equal
$1$ and $r_1^i (i=1,2), $ or $r_1$ and $r_1$, respectively. Then
$n=r_1$, $t=1$ and $r_1=3$. Thus the group ${\rm PSL}_{3}(q)$ and
$3\parallel (q-1)$, that is ${\rm PSL}_3(q)$ and $q\equiv 4,
7(\bmod~9)$.
\medskip

\emph{Case $2$: $L\cong {\rm PSU}_n(q^2)$ with $n\geq 3$.} Then
every maximal torus $T$ of ${\rm PSU}_n(q^2)$ has the order
$\frac{1}{(n,q+1)(q+1)}(q^{n_1}-(-1)^{n_1})(q^{n_2}-(-1)^{n_2})\cdots
(q^{n_k}-(-1)^{n_k})$ for appropriate partition $n=n_1+n_2+\cdots
+n_k$. Now we first assume that $m\neq 2$. If $m$ is an odd number,
then there exists a maximal torus whose order is
$\frac{1}{(n,q+1)(q+1)}(q^{2m}-1)^{[\frac{n}{2m}]}$ by the partition
$n=2m+2m+\cdots+2m+l$ where $0\le l<2m$. Furthermore, the $r_m$-part
of this order is one of the largest order of abelian $r_m$-subgroup.
Since $e_{{\rm PSU}_n(q)}(m)=[\frac{n}{2m}]$, it follows that Sylow
$r_m$-subgroups are abelian if and only if $\widetilde{e}_{L}(m)=0$,
that is the condition $e_{L}(mr_m)=0$. If $4|m$, we choose the
partition $n=m+\cdots+m+l$ and $0\le l<m$, then the order of
corresponding torus is
$\frac{1}{(n,q+1)(q+1)}(q^{m}-1)^{[\frac{n}{m}]}$. If $2<m\equiv
2(\bmod~4)$, then there is the order of a torus
$\frac{1}{(n,q+1)(q+1)}(q^{\frac{m}{2}}+1)^{[\frac{2n}{m}]}$ by the
partition $n=\frac{m}{2}+\cdots +\frac{m}{2}+l$ and $0\le
l<\frac{m}{2}$. So the largest order of abelian $r_m$-subgroup is
$e_{L}(m)$, and then the Sylow $r_m$-subgroups are abelian if and
only if
  $e_{L}(mr_m)=0$.

  Next we let $m=2.$  If $r_2\nmid n$, then there exists a
maximal torus has order $\frac{1}{(n,q+1)}(q+1)^{n-1}$ by the
partition $n=1+1+\cdots +1$. Since $e_{L}(2)=n-1$, it follows that
Sylow $r_2$-subgroups are abelian if and only if
$\widetilde{e}_{L}(2)=0$, that is $e_{L}(2r_2)=0$. When $r_2|n$,
since Sylow $r_2$-subgroup of ${\rm SU}_n(q^2)$ is one of ${\rm
SL}_n(q^2)$,   the Sylow $r_2$-subgroup of ${\rm PSU}_n(q^2)$ is one
of the group $\overline{M}:{\rm Sym}(n)$, where $\overline{M} $ is a
abelian group of order $\frac{(q^2-1)^{n-1}}{(n,q+1)}$. Assume that
$r_2^t\parallel (q+1)$. Then Sylow $r_1$-subgroups are abelian if
and only if $n=r_2=3$ and $t=1$, that is the group ${\rm
PSU}_3(q^2)$ and $q\equiv 2, 5(\bmod~9)$.

 \medskip

\emph{Case $3$: $L\cong {\rm PSp}_{2n}(q)(n\ge 3), {\rm
P\Omega}_{2n+1}(q)(n\ge 2)$.} Every maximal torus has  order
$\frac{1}{(2,q-1)}(q^{n_1}-1)(q^{n_2}-1)\cdots
(q^{n_k}-1)(q^{l_1}+1)(q^{l_2}+1)\cdots (q^{l_t}+1)$ for appropriate
partition $n=n_1+n_2+\cdots +n_k+l_1+l_2+\cdots +l_t$. Now we choose
the partition $n=m+\cdots+m+l$ if $m$ is odd,  and
$n=\frac{m}{2}+\cdots+\frac{m}{2}$ if $m$ even, where $0\le l<
\frac{m}{(2,m)}$. So   the largest order of abelian $r_m$-subgroup
is is a divisor of $\Phi_m(q)^{e_{L}(m)}$, and then the Sylow
$r_m$-subgroups are abelian if and only if
  $e_{L}(mr_m)=0$.\medskip

\emph{Case $4$: $L\cong {\rm P\Omega}_{2n}^{\epsilon}(q),
\epsilon\in\{+,-\}.$} Then every maximal torus has the order
$\frac{1}{(2,q^n-\epsilon 1)}(q^{n_1}-1)(q^{n_2}-1)\cdots
(q^{n_k}-1)(q^{l_1}+1)(q^{l_2}+1)\cdots (q^{l_t}+1)$ for appropriate
partition $n=n_1+n_2+\cdots +n_k+l_1+l_2+\cdots +l_t$ of $n$, where
$t$ is even if $\epsilon=+$ and $t$ is odd if $\epsilon=-$. So
whenever $r_m\neq 2$, it see that Sylow $r_m$-subgroup is abelian if
and only if $e_L(mr_m)=0$ by the previous method.\medskip

\emph{Case $5$: $L\cong $ $^2{\rm B}_2(q), q=2^{2n+1}.$} Since there
exist cyclic subgroups of order $q-1$, $2^{2n+1}+2^{n+1}+1$ and
$2^{2n+1}-2^{n+1}+1$ (ref. \cite[Theorem 4.2]{WR}), we have Sylow
$r$-subgroup are cyclic except $r=2$. On the other hand,
$e_L(mr_m)=0$ by Lemma 3.6.\medskip

\emph{Case $6$: $L\cong$ $^3{\rm D}_4(q).$} By Theorem 4.1 of
\cite{WR}, the simple group $^3{\rm D}_4(q)$ has  subgroups ${\rm
SL}_2(q^3)$, ${\rm SL}_3(q)$, ${\rm SU}_3(q^2)$, $(C_{q^2+q+1}\times
C_{q^2+q+1}):{\rm SL}_2(3)$, $(C_{q^2-q+1}\times C_{q^2-q+1}):{\rm
SL}_2(3)$ and $C_{q^4-q^2+1}$. Let $r_m^t||(q^m-1)$. If $m=1$ and
$r_1\neq 3$, then there is abelian subgroup $C_{q-1}\times C_{q-1}$
of ${\rm SL}_3(q)$, and so Sylow $r_1$-subgroup of $^3{\rm D}_4(q)$
is abelian. When $r_1=3$, the Sylow $r_1$-subgroup is not abelian.
If $m=2$ and $r_2\neq 3$, then there is a subgroup $C_{q+1}\times
C_{q+1}$, and so Sylow $r_2$-subgroups are abelian. If $r_2=3$, then
Sylow 3-subgroup is not abelian by above listed subgroups. If
$m=3,6,12$, then Sylow subgroups are $C_{r_3^t}\times C_{r_3^t}$,
$C_{r_6^t}\times C_{r_6^t}$ and $C_{r_{12}^t}$, respectively. These
are all abelian. Thus the Sylow $r_m$-subgroup is abelian if and
only if $e_L(mr_m)=0$.\medskip

\emph{Case $7$: $L\cong {\rm G}_2(q).$} By the section 4.2.6 of
\cite{WR} there exist subgroups ${\rm GL}_2(q)$, ${\rm SL}_3(q)$ and
${\rm SU}_3(q)$ of ${\rm G}_2(q)$, and then if $r\neq r_1(=3)$, then
Sylow $r$-subgroup is abelian. Now if $r_1=3$, then Sylow 3-subgroup
of ${\rm G}_2(q)$ is not abelian. So Sylow $r_m$-subgroup is abelian
if and only if $e_L(mr_m)=0$.\medskip

\emph{Case $8$: $L\cong$ $^2{\rm G}_2(q)$ and $q=3^{2n+1}\geq 27.$}
In terms of Table 3, we know that $e_L(mr_m)=0$ for any $m$.
Moreover, Sylow $r_m$-subgroups are cyclic by Theorem 4.3 of
\cite{WR}, and so Sylow $r_m$-subgroup is abelian if and only if
$e_L(mr_m)=0$.\medskip

\emph{Case $9$: $L\cong {\rm F}_4(q).$} Since $^3{\rm D}_4(q)$ is a
subgroup of ${\rm F}_4(q)$ (ref.
 \cite[Section 4.8]{WR}), we have Sylow $r_6$ and
$r_{12}$-subgroup of ${\rm F}_4(q)$ are ones of $^3{\rm D}_4(q)$
(see Table 1 and Table 2), which are abelian by the above Case 6.
Furthermore, ${\rm F}_4(q)$ has a subgroup ${C_2^{2}}. {\rm
P\Omega}_8^+(q)$ (ref.
 \cite[Section 4.8]{WR}), then Sylow $r_4$ and $r_8$-subgroups of
${\rm P\Omega}_8^+(q)$ are ones of ${\rm F}_4(q)$, which are abelian
by the above Case 4. Since ${\rm F}_4(q)$ has the maximal torus
$C_{q\pm 1}^4$ and $C_{q^2-q+1}^2$ (also ref. \cite[Section
4.8]{WR}), Sylow $r_m$-subgroup is abelian if and only if
$e_L(mr_m)=0$ for $m=1,2,3$.\medskip

\emph{Case $10$: $L\cong$ $^2{\rm F}_4(q)$ and $q=2^{2n+1}\ge 8$.}
Since there exists a subgroup ${\rm SU}_3(q)$(ref.  \cite[Theorem
4.4]{WR}) and Sylow $r_2$, $r_6$-subgroups are ones of $^2{\rm
F}_4(q)$ (see Table 1 and Table 2), by the discussion of Case 2, we
have these Sylow subgroups are abelian if and only if $e_L(mr_m)=0$.
Moreover, $^2{\rm F}_4(q)$ has subgroups ${\rm Sz}(q)\times
C_{q-1}$, $C_{2^{2n+1}\pm 2^{n+1}+1}^2$ and $C_{q^2+q+1\pm
2^{n+1}(q+1)}$ (also ref. \cite[Theorem 4.4]{WR}), let $r_m
\parallel q^m-1$, and then the Sylow $r_i$-subgroup is $C_{r_i^t}^2$
for $i=1, 4$, and Sylow $r_{12}$-subgroup is cyclic. So Sylow
$r_m$-subgroup is abelian if and only if $e_L(mr_m)=0$.
\medskip

\emph{Case $11$: $L\cong$${\rm E}_6(q).$} Since ${\rm F}_4(q)$ is a
subgroup of ${\rm E}_6(q)$ (ref. \cite[Section 4.6.4]{WR}), we have
Sylow $r_i$-subgroups of ${\rm E}_6(q)$ are ones of ${\rm F}_4(q)$
for $i=2,4,6,8,12$ (see Table 1 and Table 2). By the proof of Case
9, it follows that the Sylow $r_m$-subgroup is abelian if and only
if $e_L(mr_m)=0$ for $m=2,4,6,8,12$. Also the group $E_6(q)$ has
subgroups ${\rm PSL}_2(q)\times {\rm PSL}_5(q)\times C_{q-1}$ if
$r_1\neq 3,5$ and ${\rm PSL}_3(q)^3$ (also ref. \cite[Section
4.6.4]{WR}). So in this case there exist abelian subgroups
$C_{r_1^t}^6$, $C_{r_3^t}^3$ and $C_{r_5^t}$ provided $r_i^t
\parallel q^i-1$ for $i=1,3,5$, but the order of the maximal torus
whose $r_1$,$r_3$ and $r_5$-part are the largest are $(q-1)^6$,
$(q^2+q+1)^3$ and $q^5-1$, respectively. Then the Sylow $r_1(\neq
3,5)$ and $r_3$-subgroup is abelian if and only if $e_L(r_1)=0$.
Also there exists a cyclic subgroup of order $q^6+q^3+1$, so Sylow
$r_9$-subgroups are cyclic. If $r_1=5$, since ${\rm E}_6(q)$ has the
subgroup ${\rm P\Omega}_{10}^+(q)$ whose Sylow $5$-subgroups are
non-abelian(see Case 4), then the Sylow $r_1(=5)$-subgroup is
abelian if and only if $e_L(r_1)=0$. Next we consider $r_1=3$, and
use the same notation of the section 4.6.1 of \cite{WR} to denote by
${\rm GE}_6(q)$
 the group of matrices which multiply the determinant by a scalar,
 so that ${\rm GE}_6(q)\cong C_3.(C_{\frac{q-1}{3}}\times {\rm E}_6(q)).C_3$.
 Then Sylow 3-subgroup of ${\rm E}_6(q)$ is not abelian.
\medskip

\emph{Case $12$: $L\cong$ $^2{\rm E}_6(q).$} Since ${\rm F}_4(q)$ is
a subgroup of $^2{\rm E}_6(q)$ (ref.  \cite[Section 4.11]{WR}), we
have Sylow $r_i$-subgroups of $^2{\rm E}_6(q)$ is ones of ${\rm
F}_4(q)$ for $i=1,3,4,8,12$. Similarly, by the proof of Case 9, it
follows that the Sylow $r_m$-subgroup is abelian if and only if
$e_L(mr_m)=0$ for $m=1,3,4,8,12$. It is known that ${\rm PSU}_6(q)<$
$^2{\rm E}_6(q)$, we have Sylow $r_{10}$-subgroups are cyclic.
Obviously, Sylow $r_{18}$-subgroups are also cyclic. Since there is
a maximal torus $C_{q^2-q+1}^3$, we have Sylow $r_{6}$-subgroups are
abelian. Next if $r_2\neq 3$, then the order of a maximal torus
whose $r_2$-part is the largest is $(q+1)^6$, and so that Sylow
$r_2$-subgroups are not abelian. If $r_2=3$, similar to the case
${\rm E}_6(q)$, then Sylow $3$-subgroup is not abelian.\medskip

\emph{Case $13$: $L\cong$${\rm E}_7(q).$} Since ${\rm E}_6(q)$ is a
subgroup of ${\rm E}_7(q)$ (ref.  \cite[Section 4.12]{WR}), Sylow
$r_m$-subgroups of ${\rm E}_6(q)$ are ones of ${\rm E}_7(q)$ for $m=
3,4,5,8,9,12,18$  (see Table 1 and Table 2). So by Case 11, it
follows that the Sylow $r_m$-subgroup is abelian if and only if
$e_L(mr_m)=0$ for $m=3,4,5,8,9,12,18$. Also there is a section of
subgroups of ${\rm E}_7(q)$ which is isomorphic to $^2{\rm E}_6(q)$
(also ref.  \cite[Section 4.12]{WR}), but Sylow $r_m$-subgroups for
$m=6,10$ of $^2{\rm E}_6(q)$ are ones of ${\rm E}_7(q)$, and then
the Sylow $r_m$-subgroup is abelian if and only if $e_L(mr_m)=0$ for
$m=6,10$. Since there exist sections ${\rm PSU}_8(q)$ and ${\rm
PSL}_7(q)$ of subgroups of ${\rm E}_7(q)$, so the Sylow
$r_m$-subgroup is abelian if and only if $e_L(mr_m)=0$ for $m=7,14$.
Next we consider $m=1,2$. Since $r_1\ne 2$, by the structure of
maximal torus we have the largest $r_1$ and $r_2$-part of abelian
subgroups are ones of $(q-1)^7$ and $(q+1)^7$, and so Sylow $r_1$
and $r_2$-subgroups are abelian if and only if
$e_L(mr_m)=0$.\medskip

\emph{Case $14$: $L\cong$${\rm E}_8(q).$} Since ${\rm E}_7(q)$ and
${\rm P\Omega}_{16}^+(q)$ are subgroups of ${\rm E}_8(q)$ (ref.
\cite[Section 4.12]{WR}), the Sylow $r_m$-subgroup of ${\rm E}_8(q)$
is abelian if and only if $e_L(mr_m)=0$ for $m=7,8,9,14,18$. Also
$^3{\rm D}_4(q)\times $ $^3{\rm D}_4(q)$ is a subgroup of ${\rm
E}_8(q)$ (also ref. \cite[Section 4.12]{WR}), then the result is
true for $m=3,6,12$. By the structure of maximal torus we have Sylow
$r_m$-subgroups  for $m=15,20,24,30$ are cyclic, and Sylow $r_5$,
$r_{10}$-subgroups  are ones of $C_{q^4+q^3+q^2+q+1}^2$ and
$C_{q^4-q^3+q^2-q+1}^2$, respectively. So Sylow $r_5$,
$r_{10}$-subgroups of ${\rm E}_8(q)$ are ablelian. For remaining
cases of $m=1,2,4$, we use the maximal torus.  Since a maximal torus
whose $r_m$-part is the largest for $m=1,2,4$ are the groups
$C_{q\pm 1}^4$ and $C_{q^2+1}^4$, so Sylow $r_m$-subgroups are
abelian if and only if $e_L(mr_m)=0$.\hfill$\blacksquare$\\

From the above proof, we can get the following corollary.\\

\noindent \textbf{Corollary 3.9.} Let $L$ be a simple group of Lie
type over the Galois field ${\rm GF}$$(q)$ and $r_m$ an odd
primitive prime of $q^m-1$ and $r_m^t\parallel (q^m-1)$. Suppose
that $L$ has an abelian Sylow $r_m$-subgroup $R$. Then $R\cong
C_{r_m^t}^{e_L(m)}$ except the following cases:

(1) $r_1=3$,  $L\cong {\rm PSL}_3(q)$ and $q\equiv 4,7(\bmod~9)$;

(2) $r_2=3$,  $L\cong {\rm PSU}_3(q^2)$ and $q\equiv 2,5(\bmod~9)$. \\

Next we discuss  simple groups of Lie type which have an elementary
abelian Sylow subgroup. In terms of Corollary 3.9, it is sufficient
to satisfy with the condition $r_m\parallel (q^m-1)$. First, we give
a
lemma.\\

\noindent \textbf{Lemma 3.10.} Let $r_m$ be a primitive prime of
$q^m-1$ and $e$ a primitive root of $r_m$.  Let $e_i\equiv
e^{\frac{(r_m-1)i}{m}}(\bmod~r_m)$, where $(i,m)=1$, $1\leq i\leq m$
and $1\leq e_i\leq r_m$.  Then

(1) $r_1\parallel (q-1) $ if and only if $q\equiv
 1+r_1, 1+2r_1, \cdots, 1+(r_1-1)r_1(\bmod~r_1^2)$,

 (2) for $m\ge 2$,  $r_m\parallel (q^m-1) $
if and only if $q\equiv e_i, e_i+r_m, e_i+2r_m, \cdots,
e_i+(r_m-1)r_m(\bmod~r_m^2)$.  \medskip

\textbf{Proof.} Let $q\equiv x(\bmod~r_m^2)$ so that $1\leq x\leq
r_m^2$. Then we have $q^m-1\equiv x^m-1(\bmod~r_m^2)$. Hence
$r_m\parallel (q^m-1)$ if and only if the order of the element $x$
in the cyclic group $C_{r_m-1}\cong \langle e\rangle$ is $m$. But
the set of elements $x$ of order $m$ in $\langle e\rangle$ is
$\{e_i|e_i\equiv e^{\frac{(r_m-1)i}{m}}(\bmod~r_m), (i,m)=1, 1\leq
i\leq m$ and $1\leq e_i\leq r_m\}$, and then $q\equiv e_i, e_i+r_m,
e_i+2r_m, \cdots, e_i+(r_m-1)r_m(\bmod~r_m^2)$ if $m\geq 2$. For the
case of  $m=1$, clearly, $e_1=1$, and if $q\equiv 1(\bmod~r_1^2)$,
then it contradicts $r_1\parallel (q-1)$. So $q\equiv
 1+r_m, 1+2r_m, \cdots,
 1+(r_1-1)r_1(\bmod~r_1^2)$.\hfill$\blacksquare$\\

\noindent \textbf{Theorem 3.11.} Let $L$ be a simple group of Lie
type over the Galois field $\rm GF$$(q)$,  $r_m$ an odd primitive
prime of $q^m-1$ and $e$ a primitive root of $r_m$.  Set $e_i\equiv
e^{\frac{(r_m-1)i}{m}}(\bmod~r_m)$, where $(i,m)=1$, $1\leq i\leq m$
and $1\leq e_i\leq r_m$. We have the following conditions:

(1) if $(L,r_1)\neq ({\rm PSL}_3(q),3)$, then $L$ has an elementary
abelian Sylow $r_1$-subgroup if and only if   $q\equiv
 1+r_1, 1+2r_1, \cdots, 1+(r_1-1)r_1(\bmod~r_1^2)$, and

 (2) for $m\geq 2$, if $(L, r_m)\neq ({\rm PSU}_3(q^2),3)$, then $L$ has an elementary abelian Sylow $r_m$-subgroup if and only
if $q\equiv e_i, e_i+r_m, e_i+2r_m, \cdots,
e_i+(r_m-1)r_m(\bmod~r_m^2)$.\\

Finally, we give an example of classification on simple groups
having
a Sylow subgroup $C_5\times C_5$.\\

\noindent \textbf{Example 3.12.} Let $L$ be a simple group of Lie
type over finite field ${\rm GF}(q)$. Suppose that $S$ has a Sylow
subgroup $C_5\times C_5$. Then $L$ is one of following groups:

 (1)   ${\rm PSL}_2(25)$,

 (2)  ${\rm PSL}_3(q)$, ${\rm PSU}_4(q^2),$ ${\rm PSp}_4(q),$ ${\rm P \Omega}_5(q), $
$^3{\rm D}_4(q)$, ${\rm G}_2(q)$, where $q\equiv 6,11,$
$16,21$$(\bmod~25)$,

 (3) ${\rm PSL}_4(q),$ ${\rm PSL}_5(q)$, ${\rm PSU}_3(q^2),$ ${\rm P \Omega}_5(q)$,
$^3{\rm D}_4(q)$, ${\rm G}_2(q)$,  where $q\equiv4,9,14,$
$19,24$$(\bmod~25)$,

  (4) ${\rm PSL}_8(q),$ ${\rm PSL}_9(q),$ ${\rm PSL}_{10}(q),$ ${\rm PSL}_{11}(q),$
${\rm PSU}_8(q^2),$ ${\rm PSU}_9(q^2),$ ${\rm PSU}_{10}(q^2),$ ${\rm
PSU}_{11}(q^2),$ $P \Omega^+_8(q),$ $ {\rm P \Omega}^+_{10}(q),$
${\rm P \Omega}^+_{12}(q),$ ${\rm PSp}_8(q)$, ${\rm PSp}_{10}(q)$,
${\rm P\Omega}_9(q),$ ${\rm P\Omega}_{11}(q),$ ${\rm P
\Omega}^-_8(q)$, ${\rm P \Omega}^-_{10}(q),$ ${\rm P
\Omega}^-_{12}(q)$, ${\rm F}_4(q)$, $^2{\rm F}_4(q)(q=2^{2n+1}\geq
8)$, ${\rm E}_6(q)$, $^2{\rm E}_6(q)$,
${\rm E}_7(q)$, where $q\equiv 2,3,7,8,12,13,17,18,22,23(\bmod~25)$.\\

\section{Sporadic Simple Groups }

Next we give  Table 4 to list the information of Sylow subgroups of
sporadic simple groups. Note that we denote by ``+" the simple group
$S$ having abelian Sylow $r$-subgroup if ``+" in the position of
$S$-row and $r$-column, the blank otherwise. By the Atlas of finite
groups \cite{CN}, we know that if $r\ge 17$ and $r\big| |S|$, then
$r\parallel |S|$, and so Sylow $r$-subgroups are cyclic. We omit the
prime divisors which are greater than or equal to 17 in Table 4.

\begin{small}\begin{center}Table 4: {\bf~The abelian Sylow $r$-subgroup of
Sporadic Simple Groups}~~~~\end{center}
~~~~~~~~~~~\begin{tabular}{|l|l|l|l|l|l|} \hline  &3 &5
&7&11&13\\\hline
 $M_{11}$ &+&+&+&&\\\hline
  $M_{12}$ & &+&+&&\\\hline
  $J_1$ & +&+&+&+&\\\hline
  $M_{22} $  & + & +&+ &+ &\\\hline
   $J_2 $  &  & +&+ & &\\\hline
    $M_{23} $  & + & +&+ & +&\\\hline
     $HS $  &  +& &+ &+ &\\\hline
      $ J_3$  &  &+ & & &\\\hline
       $ M_{24}$  &  &+ & +&+
        &\\\hline
        $ M^cL$  &  & &+ & +&\\\hline
   $He $  &  & + & & &\\\hline
    $ Ru$  &  & & & +&+\\\hline
     $Suz $  &  & &+
      &+ &+\\\hline
\end{tabular}
~\begin{tabular}{|l|l|l|l|l|l|} \hline  &3 &5 &7&11&13\\\hline
      $ O'N$  &  &+ & &+ &\\\hline
       $ Co_3$  &  & & +&+ &\\\hline
        $Co_2 $  &  & &+ &+ &\\\hline
   $ Fi_{22}$  &  &+ &+ &+ &\\\hline
    $ HN$  &  & &+
    &+ &\\\hline
     $Ly $  &  & &+ &+ &\\\hline
      $ Th$  &  & &+ & &+\\\hline
       $ Fi_{23}$  &  &+ &+ &+ &\\\hline
        $ Co_1$  &  & & +& +&+\\\hline
   $J_4 $  &  &+ &+ & &\\\hline
    $ Fi'_{24}$  &  &+ & &+ &+\\\hline
     $B $  &  & & +& +&+\\\hline
      $ M$  &  & & &+ &\\\hline
\end{tabular}\end{small}\\\medskip

Moreover, we can know that the abelian Sylow subgroups must be an
elementary abelian group for 26 sporadic simple groups.
\\\medskip

\textbf{Acknowledgements.} Authors would like to thank two anonymous
referees for many helpful comments and suggestions that directly
lead to the improvement of the original manuscript.

 {\small


\begin{thebibliography}{36}





\bibitem{A} Artin E.,   The Orders of the Linear Groups, Comm. Pure and Applied
Math. Soc., 1959, 88-99.

\bibitem{C1}  Carter R.W., Centralizer of Semisimple Elements in the
Finite Classical Groups, Proc. London Math. Soc., 42(3), 1981, 1-41.

\bibitem{C2}  Carter R.W., Conjugacy Classes in the Weyl group, Compositio Mathematica, 25(1), 1972, 1-59.

\bibitem{CN}  Conway J.H.,  Curtis R.T.,  Norton  S.P., etc.,
Atlas of Finite Groups, Oxford: Clarendon Press, 1985.

\bibitem{F}  Feit W., On Large Zsigmondy Primes, Proc. Amer. Math.
Soc., 102(1), 1988, 29-36.

\bibitem{KLST}  Kimmerle W.,  Lyons R.,  Sandling R.,  Teague D.N.,
Composition Factors from the Group Ring and Artin's Theorem on
Orders of Simple Groups, Proc. London Math. Soc, 60(3), 1990,
89-122.


\bibitem{M}  Malle G.,   Moreto A.,  Navarro G., Element Orders and Sylow
Structure, Mathematische Zeitschrift, 252(1), 2006, 223-230.

\bibitem{VV}  Mal'tsev A.I.,  Commutative Subalgebras of Semisimple Lie Algebras, Izv. Akad. Nauk SSSR, Ser.
Mat., 9(4),1945, 291-300.


\bibitem{R}   Ribenboin P., The Book of Prime Number Records, Second
Edition, Springer-Verlag, New York, 1989.



\bibitem{V}  Vdovin E.P., Large Abelian Unipotent Subgroups of Finite Chevalley
Groups, Algebra and Logic,  40 (5), 2001, 292-305.



\bibitem{Wa}  Walter J.H., The Characterzation of Finite Groups with
Abelian Sylow 2-Subgroup, Ann. Math., 89, 1969, 405-514.

\bibitem{Weir}  Weir A.J., Sylow $p$-subgroups of the Classical Groups over Finite Fields with Characteristic Prime to $p$,  Proc. Amer. Math. Soc., 6, 1955, 529-533.


\bibitem{WR} Wilson R.A., The Finite Simple Groups, Graduate Texts in Mathematics 251, Springer-Verlag, 2009.

\bibitem{Z}  Zsigmondy K., Zur Theorie der Potenzreste, Monatsch.
Math. Phys., 3(1892), 265-284.




\end{thebibliography}
\end{document}